\documentclass[11pt]{amsart}

\newtheorem{theorem}{Theorem}[section]

\theoremstyle{definition}
\newtheorem{definition}[theorem]{Definition}

\theoremstyle{remark}

\numberwithin{equation}{section}

\begin{document}

\newcommand{\spacing}[1]{\renewcommand{\baselinestretch}{#1}\large\normalsize}
\spacing{1.14}

\title{On some hypercomplex 4-dimensional Lie groups of constant scalar curvature}

\author {H. R. Salimi Moghaddam}

\address{Department of Mathematics, Faculty of  Sciences, University of Isfahan, Isfahan,81746-73441-Iran.} \email{salimi.moghaddam@gmail.com and hr.salimi@sci.ui.ac.ir}



\begin{abstract}
In this paper we study sectional curvature of invariant
hyper-Hermitian metrics on simply connected $4$-dimensional real
Lie groups admitting invariant hypercomplex structure. We give the
Levi-Civita connections and explicit formulas for computing
sectional curvatures of these metrics and show that all these
spaces have constant scalar curvature. We also show that they are
flat or they have only non-negative or non-positive sectional
curvature.
\end{abstract}

\maketitle

\section{\textbf{Introduction}}\label{Intro}
Hypercomplex structures on $4n$-dimensional manifolds are of
interesting structures in mathematics which have many applications
in theoretical physics. For example these structures appear in
supersymmetric sigma model \cite{Po}. Gibbons, Papadopoulos and
Stelle showed that the geometry of the moduli space of a class of
black holes in five dimensional is hyper-K$\ddot{a}$hler
\cite{GiPaSt}. Therefore it is important to study hypercomplex
spaces. In this paper we consider invariant hyper-Hermitian
metrics on simply connected $4$-dimensional real Lie groups
admitting invariant hypercomplex structure. These spaces
classified by M. L. Barberis (see \cite{Ba1}). One of important
quantities which associate to Riemannian manifold is sectional
curvature. In this article after obtaining the Levi-Civita
connections of these Riemannian spaces we give explicit formulas
for computing sectional curvatures of these manifolds. Then by
using these formulas we show that these spaces have constant
scalar curvature.

\section{\textbf{Preliminaries}}\label{Prelim}

\begin{definition}
Suppose that $M$ is a $4n$-dimensional manifold. Also let $J_i,
i=1,2,3,$ be three fiberwise endomorphism of $TM$ such that
\begin{eqnarray}
  J_1J_2&=& -J_2J_1=J_3, \label{JJ} \\
  J_i^2 &=& -Id_{TM},  \ \ \ \ \ \ \ i=1,2,3, \\
  N_i &=& 0, \ \ \ \ \ \ \ i=1,2,3,
\end{eqnarray}
where $N_i$ is the Nijenhuis tensor (torsion) corresponding to
$J_i$ defined as follows:
\begin{eqnarray}
   N_i(X,Y)=[J_iX,J_iY]-[X,Y]-J_i([X,J_iY]+[J_iX,Y]),
\end{eqnarray}
for all vector fields $X,Y$ on $M$. The family
$\mathcal{H}=\{J_i\}_{i=1,2,3}$ is a hypercomplex structure on
$M$.
\end{definition}
In fact a hypercomplex structure on a $4n$-dimensional manifold
$M$ is a family $\mathcal{H}=\{J_i\}_{i=1,2,3}$ of complex
structures on $M$ satisfying in the relation \ref{JJ} (since an
almost complex structure is a complex structure if  and only if it
has no torsion, see \cite{KoNo} page 124.).
\begin{definition}
A Riemannian metric $g$ on a hypercomplex manifold
$(M,\mathcal{H})$ is called hyper-Hermitian if
$g(J_iX,J_iY)=g(X,Y)$, for all vector fields $X,Y$ on $M$ and
$i=1,2,3$.
\end{definition}
\begin{definition}
A hypercomplex structure $\mathcal{H}=\{J_i\}_{i=1,2,3}$ on a Lie
group $G$ is said to be (left) invariant if for any $a\in G$,
\begin{eqnarray}
  J_i=Tl_a\circ J_i\circ Tl_{a^{-1}},
\end{eqnarray}
where $Tl_a$ is the differential function of the left translation
$l_a$.
\end{definition}
From now we suppose that $G$ is a simply connected $4$-dimensional
real Lie group.\\
Suppose that $g$ is a left invariant Riemannian metric on a Lie
group $G$ with Lie algebra $\frak{g}$, then the Levi-Civita
connection of $g$ is defined by the following relation
\begin{eqnarray}\label{nabla}
  2g(\nabla_UV,W)=g([U,V],W)-g([V,W],U)+g([W,U],V),
\end{eqnarray}
for any $U,V,W\in\frak{g}$, where $<,>$ is the inner product
induced by $g$ on $\frak{g}$.


\section{\textbf{Sectional and Scalar Curvatures}}
In this section we compute the Levi-Civita connections and the
sectional curvatures of left invariant hyper-Hermitian metrics on
(left) invariant hypercomplex $4$-dimensional simply connected Lie
groups. Then we compute the scalar curvatures of these spaces and
show that these spaces are of constant scalar curvature. We
mention that these spaces are complete because they are
homogeneous Riemannian manifolds. M. L. Barberis has classified
these Lie
groups in \cite{Ba1}.\\
Let $G$ be a Lie group as above with Lie algebra $\frak{g}$. She
has shown that $g$ is either Abelian or isomorphic to one of the
following Lie algebras:
\begin{enumerate}
    \item $[Y,Z]=W$, $[Z,W]=Y$, $[W,Y]=Z$, $X$ central,
    \item $[X,Z]=X$, $[Y,Z]=Y$, $[X,W]=Y$, $[Y,W]=-X$,
    \item $[X,Y]=Y$, $[X,Z]=Z$, $[X,W]=W$,
    \item $[X,Y]=Y$, $[X,Z]=\frac{1}{2}Z$, $[X,W]=\frac{1}{2}W$,
    $[Z,W]=\frac{1}{2}Y$,
\end{enumerate}
where $\{X,Y,Z,W\}$ is an orthonormal basis.\\
All above Lie groups are diffeomorphic to $\Bbb{R}^4$ unless case
(1) which is diffeomorphic to $\Bbb{R}\times S^3$ (see \cite{Ba1}
and \cite{Ba2}.).\\
Now we compute the Levi-Civita connection and sectional curvature
of any case separately.\\
In Abelian case obviously $(G,g)$ is flat therefore we consider
the other cases.\\

\textbf{Case 1.} A direct computation, by using formula
\ref{nabla}, shows that for Levi-Civita connection of this case we
have:
\begin{eqnarray}
  &&\nabla_XX = 0 \ \ , \ \ \nabla_XY = 0, \ \ \nabla_XZ = 0, \ \ \nabla_XW = 0, \nonumber\\
  &&\nabla_YX = 0 \ \ , \ \ \nabla_YY = 0, \ \ \nabla_YZ = \frac{1}{2}W, \ \ \nabla_YW = -\frac{1}{2}Z,\\
  &&\nabla_ZX = 0 \ \ , \ \ \nabla_ZY = -\frac{1}{2}W, \ \ \nabla_ZZ = 0, \ \ \nabla_ZW = \frac{1}{2}Y,\nonumber\\
  &&\nabla_WX = 0 \ \ , \ \ \nabla_WY = \frac{1}{2}Z \ \ , \ \ \nabla_WZ = -\frac{1}{2}Y, \ \ \nabla_WW = 0\nonumber.
\end{eqnarray}
Now by using Levi-Civita connection for curvature tensor we have:
\begin{eqnarray}
  &&R(Y,Z)Y = -R(Z,W)W = -\frac{1}{4}Z,\nonumber\\
  &&R(Y,W)W = R(Y,Z)Z = \frac{1}{4}Y,\\
  &&R(Z,W)Z = R(Y,W)Y = -\frac{1}{4}W\nonumber,
\end{eqnarray}
and in other cases $R=0$.

Let $U=aX+bY+cZ+dW$ and
$V=\tilde{a}X+\tilde{b}Y+\tilde{c}Z+\tilde{d}W$ be two arbitrary
vectors in $\frak{g}$ then we have:
\begin{eqnarray}
  R(V,U)U=-\frac{1}{4}\{(b\tilde{c}-c\tilde{b})(cY-bZ)+(b\tilde{d}-d\tilde{b})(dY-bW)+(c\tilde{d}-d\tilde{c})(dZ-cW)\}.
\end{eqnarray}
Now assume that $\{U,V\}$ is an orthonormal set then for sectional
curvature $K(U,V)$ we have:
\begin{eqnarray}
  K(U,V) =
  \frac{1}{4}\{(b\tilde{c}-c\tilde{b})^2+(b\tilde{d}-d\tilde{b})^2+(c\tilde{d}-d\tilde{c})^2\}\geq0.
\end{eqnarray}
The last equation shows that in this case $(G,g)$ is of
non-negative sectional curvature.\\
Now consider the orthonormal basis $\{X,Y,Z,W\}$. By using above
formula for sectional curvature we can show that the scalar
curvature $S$ of this space is $S=\frac{3}{2}$.

\textbf{Case 2.} Similar to case 1 by using \ref{nabla} we have
\begin{eqnarray}
  &&\nabla_XX = -Z \ \ , \ \ \nabla_XY = 0, \ \ \nabla_XZ = X, \ \ \nabla_XW = 0, \nonumber\\
  &&\nabla_YX = 0 \ \ , \ \ \nabla_YY = -Z, \ \ \nabla_YZ = Y, \ \ \nabla_YW = 0,\\
  &&\nabla_ZX = 0 \ \ , \ \ \nabla_ZY = 0, \ \ \nabla_ZZ = 0, \ \ \nabla_ZW = 0,\nonumber\\
  &&\nabla_WX = -Y \ \ , \ \ \nabla_WY = X \ \ , \ \ \nabla_WZ = 0, \ \ \nabla_WW = 0\nonumber.
\end{eqnarray}
The curvature tensor of the above connection is as follows:
\begin{eqnarray}
  &&R(X,Y)X = -R(Y,Z)Z = Y,\nonumber\\
  &&R(X,Y)Y = R(X,Z)Z = -X,\\
  &&R(X,Z)X = R(Y,Z)Y = Z\nonumber,
\end{eqnarray}
and in other cases $R=0$. In this case for any $U$ and $V$ we
have:
\begin{eqnarray}
  R(V,U)U=-\{(a\tilde{b}-b\tilde{a})(aY-bX)+(a\tilde{c}-c\tilde{a})(aZ-cX)+(b\tilde{c}-c\tilde{b})(bZ-cY)\},
\end{eqnarray}
and for an orthonormal set $\{U,V\}$ the sectional curvature
$K(U,V)$ is
\begin{eqnarray}
  K(U,V) =
  -\{(a\tilde{b}-b\tilde{a})^2+(a\tilde{c}-c\tilde{a})^2+(b\tilde{c}-c\tilde{b})^2\}
  \leq 0,
\end{eqnarray}
which shows that in the case 2 $(G,g)$ is of non-positive
sectional curvature.\\
Now the formula of sectional curvature shows that the scalar
curvature is of the form $S=-6$.

\textbf{Case 3.} If we repeat the computations for case 3 then we
have:
\begin{eqnarray}
  &&\nabla_XX = 0 \ \ , \ \ \nabla_XY = 0, \ \ \nabla_XZ = 0, \ \ \nabla_XW = 0, \nonumber\\
  &&\nabla_YX = -Y \ \ , \ \ \nabla_YY = X, \ \ \nabla_YZ = 0, \ \ \nabla_YW = 0,\\
  &&\nabla_ZX = -Z \ \ , \ \ \nabla_ZY = 0, \ \ \nabla_ZZ = X, \ \ \nabla_ZW = 0,\nonumber\\
  &&\nabla_WX = -W \ \ , \ \ \nabla_WY = 0 \ \ , \ \ \nabla_WZ = 0, \ \ \nabla_WW =
  X\nonumber,
\end{eqnarray}
and therefore for $R$ we have
\begin{eqnarray}
  &&R(X,Y)X = -R(Y,Z)Z = -R(Y,W)W = Y,\nonumber\\
  &&R(X,Y)Y = R(X,Z)Z = R(X,W)W = -X,\\
  &&R(X,Z)X = R(Y,Z)Y = -R(Z,W)W = Z\nonumber,\\
  &&R(X,W)X = R(Y,W)Y = R(Z,W)Z = W\nonumber,
\end{eqnarray}
and in other cases $R=0$. Then for any $U$ and $V$ we have:
\begin{eqnarray}
  R(V,U)U&=&-\{(a\tilde{b}-b\tilde{a})(aY-bX)+(a\tilde{c}-c\tilde{a})(aZ-cX)+(b\tilde{c}-c\tilde{b})(bZ-cY)\nonumber\\
  &&+(a\tilde{d}-d\tilde{a})(aW-dX)+(b\tilde{d}-d\tilde{b})(bW-dY)+(c\tilde{d}-d\tilde{c})(cW-dZ)\}.
\end{eqnarray}
Hence for an orthonormal set $\{U,V\}$ the sectional curvature
$K(U,V)$ is as follows:
\begin{eqnarray}
  K(U,V)&=&-\{(a\tilde{b}-b\tilde{a})^2+(a\tilde{c}-c\tilde{a})^2+(b\tilde{c}-c\tilde{b})^2\nonumber\\
  &&+(a\tilde{d}-d\tilde{a})^2+(b\tilde{d}-d\tilde{b})^2+(c\tilde{d}-d\tilde{c})^2\} \leq 0.
\end{eqnarray}
Therefore $(G,g)$ is of non-positive sectional curvature.\\
Similar to above cases we can obtain that the scalar curvature is
of the form $S=-12$.

\textbf{Case 4.} Similar to aforementioned cases we can obtain
$\nabla$ and $R$ as follows:
\begin{eqnarray}
  &&\nabla_XX = 0 \ \ , \ \ \nabla_XY = 0, \ \ \nabla_XZ = 0, \ \ \nabla_XW = 0, \nonumber\\
  &&\nabla_YX = -Y \ \ , \ \ \nabla_YY = X, \ \ \nabla_YZ = -\frac{1}{4}W, \ \ \nabla_YW = \frac{1}{4}Z,\\
  &&\nabla_ZX = -\frac{1}{2}Z \ \ , \ \ \nabla_ZY = -\frac{1}{4}W, \ \ \nabla_ZZ = \frac{1}{2}X, \ \ \nabla_ZW = \frac{1}{4}Y,\nonumber\\
  &&\nabla_WX = -\frac{1}{2}W \ \ , \ \ \nabla_WY = \frac{1}{4}Z \ \ , \ \ \nabla_WZ = -\frac{1}{4}Y, \ \ \nabla_WW =
  \frac{1}{2}X\nonumber,
\end{eqnarray}

\begin{eqnarray}
  &&-R(X,Y)Y = -4R(X,Z)Z = -4R(X,W)W = 8R(Y,Z)W = -8R(Y,W)Z = -4R(Z,W)Y = X,\nonumber\\
  &&R(X,Y)X = -8R(X,Z)W = 8R(X,W)Z = -\frac{16}{7}R(Y,Z)Z = -\frac{16}{7}R(Y,W)W = 4R(Z,W)X =Y\nonumber,\\
  &&-4R(X,Y)W = 4R(X,Z)X = -8R(X,W)Y = \frac{16}{7}R(Y,Z)Y = 8R(Y,W)X = -\frac{16}{7}R(Z,W)W =Z\nonumber,\\
  &&4R(X,Y)Z = 8R(X,Z)Y = 4R(X,W)X = -8R(Y,Z)X = \frac{16}{7}R(Y,W)Y = \frac{16}{7}R(Z,W)Z =W\nonumber,
\end{eqnarray}
and in other cases $R=0$. For any $U$ and $V$ we have:
\begin{eqnarray}
  R(V,U)U&=&-\{(a\tilde{b}-b\tilde{a})(aY-bX-\frac{d}{4}Z+\frac{c}{4}W)
                +(a\tilde{c}-c\tilde{a})(-\frac{c}{4}X-\frac{d}{8}Y+\frac{a}{4}Z+\frac{b}{8}W)\nonumber\\
            &&+(a\tilde{d}-d\tilde{a})(-\frac{d}{4}X+\frac{c}{8}Y-\frac{b}{8}Z+\frac{a}{4}W)
                +(b\tilde{c}-c\tilde{b})(\frac{d}{8}X-\frac{7c}{16}Y+\frac{7b}{16}Z-\frac{a}{8}W)\\
            &&+(b\tilde{d}-d\tilde{b})(-\frac{c}{8}X-\frac{7d}{16}Y+\frac{a}{8}Z+\frac{7b}{16}W)
                +(c\tilde{d}-d\tilde{c})(-\frac{b}{4}X+\frac{a}{4}Y-\frac{7d}{16}Z+\frac{7c}{16}W)\}\nonumber.
\end{eqnarray}
The above equation shows that for an orthonormal set $\{U,V\}$ the
sectional curvature $K(U,V)$ is as follows:
\begin{eqnarray}
  K(U,V)&=&-\{(a\tilde{b}-b\tilde{a}+\frac{1}{4}(c\tilde{d}-d\tilde{c}))^2
               +\frac{1}{4}(a\tilde{c}-c\tilde{a}+\frac{1}{2}(b\tilde{d}-d\tilde{b}))^2\nonumber\\
          &&+\frac{1}{4}(a\tilde{d}-d\tilde{a}-\frac{1}{2}(b\tilde{c}-c\tilde{b}))^2
               +\frac{3}{8}((b\tilde{c}-c\tilde{b})^2+(b\tilde{d}-d\tilde{b})^2+(c\tilde{d}-d\tilde{c})^2)\}\leq 0.
\end{eqnarray}
Therefore $(G,g)$ is of non-positive sectional curvature.\\
In this case the scalar curvature is $S=-\frac{45}{8}$.

\large{\textbf{Acknowledgment}}\\
This work was supported by the research grant from Shahrood
university of technology.


\bibliographystyle{amsplain}

\end{document}